\documentclass[]{amsart}
\newtheorem{theorem}{Theorem}[section]
\newtheorem{thm}[theorem]{Theorem}
\newtheorem{lem}[theorem]{Lemma}
\newtheorem{prop}[theorem]{Proposition}
\newtheorem{cor}[theorem]{Corollary}

\theoremstyle{definition}

\theoremstyle{remark}
\newtheorem{rem}[theorem]{Remark}

\theoremstyle{remark}

\theoremstyle{remark}
\newtheorem{exam}[theorem]{Example}
\newcommand{\fm}{\mathfrak m}
\newcommand{\la}{\lambda}
\newcommand{\depth}{\text{depth}}

\begin{document}

\title[Fiber cone of ideals with analytic spread one]
{On the structure of the fiber cone of \\ ideals with analytic spread one}
\author{Teresa Cortadellas and Santiago Zarzuela}
\address{Departament d'{\`A}lgebra i Geometria \\ Gran Via de les Corts Catalanes 585 \\ 08007 Barcelona \\ Spain}
\email{terecortadellas@ub.edu \\ szarzuela@ub.edu}
\thanks{Both authors supported by MTM2004-01850 (Spain)}
\subjclass[2000]{Primary 13A30; Secondary 13H10, 13H15, 13A02}
\begin{abstract}
For a given a local ring $(A, \mathfrak m)$, we study the fiber
cone of ideals in $A$ with analytic spread one. In this case, the
fiber cone has a structure as a module over its Noether
normalization which is a polynomial ring in one variable over the
residue field. One may then apply the structure theorem for
modules over a principal domain to get a complete description of
the fiber cone as a module. We analyze this structure in order to
study and characterize in terms of the ideal itself the
arithmetical properties and other numerical invariants of the
fiber cone as multiplicity, reduction number or
Castelnuovo-Mumford regularity.
\end{abstract}

\maketitle{}

\section{Introduction}
\label{A}

Let $(A,\fm)$ be a Noeherian local ring and let $I$ be an ideal of
$A$. The \textit{fiber cone of $I$} (or the special fiber of the
Rees algebra $A[It]$) is the ring
$$F(I)=\bigoplus_{n\geq 0}I^n/\fm I^n =A[It]\otimes _A A/\fm .$$
Its Krull dimension is called the \textit{analytic spread of $I$}
and we will denote it by $l(I)$.

An ideal $J\subseteq I$ is called a \textit{reduction of $I$} if
there exists an integer $n$ such that $I^{n+1}=JI^n$. Phrased
otherwise, $J$ is a reduction of $I$ if
$$A[Jt] \hookrightarrow A[It]$$
is a finite morphism of graded algebras. Equivalently, it is known
that $J$ is a reduction of $I$ if and only if $I$ is integral over
$J$.

A reduction $J$ of $I$ is a minimal reduction if $J$ is minimal
with respect to inclusion among reductions of $I$. By Northcott
and Rees \cite{NoRe} minimal reductions always exist. Let $J$ be a
reduction of $I$ and assume in addition that the residue field of
$A$ is infinite. Then, $J$ is a minimal reduction of $I$ if, and
only if, $J$ is minimally generated by $l(I)$ elements if, and
only if, $J$ is generated by a family of analytically independent
elements in $I$. Therefore, given $J$ a minimal reduction of $I$,
the ring $F(J)$ is isomorphic to a polynomial ring in $l(I)$
variables over $A/\fm$ and the equalities $\fm I^n \cap J^n=\fm
J^n$ are satisfied for all $n$. That is, the graded morphism
$$F(J) \hookrightarrow F(I)$$
is a Noether normalization.

For $a\in I$, we will denote by $a^0$ the class of $a$ in $I/\fm
I$. Minimal reductions also provide homogeneous systems of
parameters of $F(I)$. Concretely, if the residue field of $A$ is
infinite, a family of elements $a_1,\dots ,a_l\in I$ is a minimal
set of generators of a minimal reduction of $I$ if and only if
$a_1^0,\dots ,a_l^0$ is a homogeneous system of parameters of
$F(I)$.

Assume now that the residue field is infinite and $l(I)=1$. If
$J=(a)$ is a minimal reduction of $I$, then $F(J)$ is isomorphic
to a polynomial ring in one variable over $A/\fm$ and $F(I)$ is a
graded finite module over $F(J)$. So we may apply the structure
theorem of finitely generated graded modules over a principal
ideal graded domain to get a set of invariants describing the
precise structure of $F(I)$ as $F(J)$-module.

Our purpose in this paper is to analyze in detail the information
provided by this set of invariants in order to study the
properties of fiber cones of dimension one. In particular, the
Cohen-Macaulay, Gorenstein or Buchsbaum properties, and other
numerical information such as Castelnuovo-Mumford regularity,
multiplicity, Hilbert funtion, reduction number or postulation
number. As we will see, although the structure of $F(I)$ as
$F(J)$-module is less rich than the structure of $F(I)$ as
$F(J)$-algebra, it suffices in this case to characterize all the
above properties in terms of the ideal itself.

The fiber cone of an ideal $I$ is one of the so called
\textit{blow up algebras of $I$} and its \emph{Proj} represents
the fiber of the maximal ideal $\fm$ by the blow up of $A$ with
center $I$. Moreover, it provides interesting information about
the ideal itself: The Hilbert function of the fiber cone describes
the minimal number of generators of the powers of $I$ and, when
the residue field is infinite, its dimension coincides with the
minimal number of generators of any minimal reduction of $I$. For
the maximal ideal itself, the fiber cone coincides with the
associated graded ring, and so in this particular situation it has
been extensively studied, the case of analytic spread one being
the tangent cones of curve singularities. But for a general ideal,
the properties of the fiber cone are much less known.
Nevertheless, in recent years some effort has been done by several
authors in order to understand its behaviour.

With respect to the arithmetical properties of the fiber cone, one
of the first known results was given by Huneke and Sally
\cite{HunSa} who proved that, if $A$ is Cohen-Macaulay, the fiber
cone of any $\fm$-primary ideal of reduction number one is
Cohen-Macaulay. This result was later extended by K. Shah
\cite{Sh1,Sh2} to equimultiple ideals of reduction number one,
giving also some conditions for the Cohen-Macaulayness of the
fiber cone of equimultiple ideals of reduction number two.
Subsequent results by Cortadellas and Zarzuela \cite{CoZa1,CoZa2},
D'Cruz, Raghavan and Verma \cite{D'CrRaVe}, and D'Cruz and Verma
\cite{D'CrVe} completed the results of Shah for more general
families of ideals. Also, the fiber cone of the defining ideal of
a monomial curve in $\mathbb{P}^3$ lying on a quadric was proven
to be Cohen-Macaulay by Morales and Simis \cite{MoSi}. This result
was later extended by P. Gimenez \cite{Gi} and Barile and Morales
\cite{BaMo} to the defining ideal of a projective monomial variety
of codimension two.

On the other hand, motivated by work of R. H\"{u}bl \cite{Hub},
H\"{u}bl and Huneke \cite{HubHun} studied the Cohen-Macaulay
property of the fiber cone of special ideals in connection with
the theory of evolutions introduced by Eisenbud and Mazur
\cite{EiMa}, which is related to A. Wiles's work on Fermat's Last
Theorem \cite{Wi}. H\"{u}bl and Swanson \cite{HubSw} have also
made some concrete computations on fiber cones in this context.
More recent work concerning the properties of the fiber cones
(multiplicity, Hilbert function, Cohen-Macaulayness,
Gorensteiness, depth...) has been done by Corso, Ghezzi, Polini
and Ulrich \cite{CoGePoUl}, Corso, Polini and Vasconcelos
\cite{CoPoVa}, T. Cortadellas \cite{Co}, D'Cruz and Puthepurakal
\cite{D'CrPu}, Heinzer and Kim \cite{HeKi}, Heinzer, Kim and
Ulrich \cite{HeKiUl}, Jayanthan and Verma \cite{JaVe1,JaVe2},
Jayanthan, Puthepurakal and Verma \cite{JaPuVe}, or D. Q. Vi{\^e}t
\cite{Vi} and others.

The case of ideals having a principal reduction has also been
considered in some detail by several authors. S. Huckaba
\cite{Huc} studied the reduction number and observed that, for a
regular ideal of analytic spread one, the reduction number doesn't
depend on the minimal reduction. And more recently, D'Anna,
Guerrieri and Heinzer \cite{D'AnGuHe1,D'AnGuHe2} have also
considered several aspects of these ideals, such as their fiber
cone, the relation type or the Ratliff-Rush closure. On the other
hand, one can also find the case of analytic spread one ideals in
induction arguments, such as the so-called Sally machine for fiber
cones, see Jayanthan and Verma \cite{JaVe2}.

Next, we briefly explain the content and structure of this paper.
In Section \ref{B} we describe the structure of $F(I)$ as
$F(J)$-module, introducing the set of invariants provided by this
structure. We relate them to several other numerical invariants of
the ideal such as reduction number or minimal number of
generators, and of the fiber cone such as multiplicity, regularity
or postulation number. Then, we give some formulas which allow to
compute this set of invariants in terms of lengths of annihilator
ideals. In particular, we prove the invariance with respect to the
chosen reduction $J$ of a distinguished subset of this set of
invariants. Section \ref{C} is devoted to the study of the
Gorenstein, Cohen-Macaulay and Buchsbaum properties of the fiber
cone; We give several characterizations of all these properties,
both in terms of the set of invariants coming form the structure
of $F(I)$ as $F(J)$-module and the corresponding lengths of
annihilator ideals introduced in the previous section. We point
out that the Buchsbaum property of the fiber cone is equivalent to
the fact that its structure as a module over $F(J)$ is independent
of the chosen minimal reduction $J$. In Section \ref{D} we give
some applications and explicit examples, which explain the results
obtained in Sections \ref{B} and \ref{C}. In particular, we get
that the fiber cone of any regular ideal with analytic spread one
and reduction number two is Buchsbaum, and give examples showing
that this is no more true for reduction number three. Finally, in
Section \ref{E} we use induction arguments to extend some of the
previous results to ideals of higher analytic spread, recovering
several known results for which we give an alternative and easier
proof.

Throughout this paper we will assume that $(A,\fm)$ is a local
Noetherian ring with an infinite residue field. For all
unexplained terminology one may use Bruns and Herzog \cite{BrHe}.

\section{The structure of $F(I)$ as $F(J)$-module}
\label{B} Let $I$ be an ideal of $A$ with analytic spread
$l(I)=l$. Let $J\subseteq I$ be a minimal reduction. Then, the
least integer $r$ such that $I^{r+1}=JI^r$ is the {\it reduction
number} of $I$ with respect to $J$ and it is denoted by $r_J(I)$.
Let  $Y_1,\dots ,Y_s$ be a minimal set of homogeneous generators
of $F(I)$ as $F(J)$-module. Then, by lifting the equality
$F(I)=\sum F(J)Y_i$ to $A[It]$ and by Nakayama's Lemma one gets
that
 $$ r_J(I)=\max \{ \deg(Y_i),\, 1\leq i \leq s \}.$$

Recall that given a finitely generated graded module $M$ over a
polynomial ring $k[x_1, \ldots ,x_n]$ over a field $k$ and a
minimal graded free resolution of $M$
$$0\rightarrow F_s \rightarrow \cdots \rightarrow F_1 \rightarrow
F_0 \rightarrow M \rightarrow 0\,,$$ the
\textit{Castelnuovo-Mumford regularity} of $M$ is the number
$$\mathrm{reg}(M):=\max \{b_i(M)-i|i=0, \ldots ,s\}\,,$$ where
$b_i(M)$ denotes the maximum of the degrees of the generators of
$F_i$.

More in general, let $S=\oplus _{n\geq 0} S_n$ be a finitely
generated standard graded algebra over a Noetherian commutative
ring $S_0$ and let $S_+=\oplus _{n>0} S_n$ be the irrelevant ideal
of $S$. Given $M=\oplus _{n\in \mathbb{Z}} M_n$ a finitely
generated graded $S$-module, let $H^i_{S_+}(M)$ be the $i$-th
graded local cohomology module of $M$ with respect to $S_+$. For
any graded $S$-module $N$ we consider
$$ \mathrm{end}(N):= \left\{
\begin{array}{ll}  \sup \{n \mid N_n\neq 0 \} & \text{ if } N\neq
0 \\ -\infty & \text{ if } N=0
\end{array} \right. $$ and denote by $a_i(M):= \mathrm{end}(H^i_{S+}(M))$.
Then, the {\it Castelnuovo-Mumford regularity} of $M$ is the
number $$\text{reg}_S(M):=\max \{ a_i(M)+i \mid i\geq 0\}.$$ It is
well known that this definition extends the classical definition
of Castelnuovo-Mumford regularity for modules over a polynomial
ring.

Observe that since $\mathrm{rad}(F(J)_+F(I))=F(I)_+$, then for
every graded $F(I)$-module $M$ one has
$H^i_{F(I)_+}(M)=H^i_{F(J)_+}(M)$ and so
$$\mathrm{reg}(F(I)):=\mathrm{reg}_{F(I)}(F(I))=\mathrm{reg}_{F(J)}(F(I)).$$

Let $S$ be a graded standard algebra over a field $k$. We shall
denote the {\it length} of a graded $S$-module $M$ by $\la (M)=
\la _{S}(M)= \la_{k}(M)=\sum \la_{k} (M_n)$. Then, $\la
(F(I))=\la_{F(I)}(F(I))=\la_{F(J)}(F(I))=\sum \la_{A/\fm}(I^n/\fm
I^n)$. Let $H(F(I),n)=\la_{A/\fm} (I^n/\fm I^n)
 =\mu (I^n)$ be the {\it Hilbert function } of $F(I)$. Then,
 $H(F(I),n)$ is of polynomial type of degree $l-1$. The unique
 polynomial $P_{F(I)}(x)\in \mathbb Q[x]$ for which
 $H(F(I),n)=P_{F(I)}(n)$ for all $n$ large enough is the \textit{Hilbert polynomial}
 of $F(I)$ and has the form
 $$P_{F(I)}(x)=\sum_{i=0}^{l-1} (-1)^{l-1-i} e_{l-1-i} {x+i\choose i}.$$
 The \textit{multiplicity} of $F(I)$ is defined as
 $$e(F(I))= \left\{ \begin{array}{lr} e_0 \quad &\text{ if } l>0 \\ \la(F(I))
 \quad &\text{ if } l=0 \end{array} \right.$$
and the {\it fiber postulation number} $\mathrm{fp}(I)$ of $I$ as
the largest integer $n$ such that $P_{F(I)}(n) \neq
H(F(I),n)=\mu(I^n)$.

 Let $H_{F(I)}(x)=\displaystyle{\sum_{n\geq 0}} \mu(I^n)x^n$ the \textit{Hilbert
series} of $F(I)$. Then
$$H_{F(I)}(x)=\frac{Q_{F(I)}(x)}{(1-x)^l}$$
for an unique $Q_{F(I)}(x) \in \mathbb Z[x,x^{-1}]$ and
$Q_{F(I)}(1)=e(F(I))$.

Assume now that $l(I)=1$ and let $J=(a)$ be a minimal reduction of
$I$. Then, $F(J)$ is isomorphic to a polynomial ring in one
variable over the residue field and so a graded principal ideal
domain. In this way, we can consider the graded decomposition of
$F(I)$ as direct sum of cyclic graded $F(J)$-modules, see also W.
V. Vasconcelos \cite[9.3]{Va},
$$F(I) \simeq \bigoplus_{i=1}^{e}F(J)(-b_i)
 \bigoplus_{j=1}^{f}(F(J)/a^{c_j}F(J))(-d_j) \quad \hfill{(\ast)}$$
where we may assume $b_1\leq \dots\leq b_e$, $d_1\leq \dots \leq
d_f$. In particular one immediately has
$$H_{F(I)}(x)=
 \frac{x^{b_1}+\cdots + x^{b_e}+(1-x^{c_1})x^{d_1}+\cdots + (1-x^{c_f})x^{d_f}}{1-x} .$$
 Moreover, in this case the Hilbert polynomial  $P_{F(I)}(x)=e(F(I))$
is a constant. As a consequence, for these ideals we have that
$F(I)$ (as a $F(J)$-module) satisfies
$$ \begin{array}{cl}
 \mu_{F(J)}(F(I))&=e+f,\\
 r_J(I)& = \max \{b_e,d_f\},\\
 \text{reg}(F(I))&= \max \{ b_e,c_j+d_j-1 \},\\
 e(F(I))&=e.
 \end{array}
 $$

Assume moreover that $I$ contains a regular element: These ideals
are usually called \textit{regular
 ideals}. One immediately gets that $a$ must be a regular element. Put
 $r:=r_J(I)$. If $n\geq r$ then
 $$ I^n/\fm I^n=a^{n-r}I^r/a^{n-r}\fm I^r \cong I^r/\fm I^r$$
 and $\mu (I^n)=\mu (I^r)$. So, the postulation number $\mathrm{fp}(I)\leq
 r-1$ and
   the Hilbert series is in this case
$$H_{F(I)}(x)=\sum _{n\geq 0} \mu (I^n) x^n=
 \frac{1+ (\mu (I)-1)x +\cdots +(\mu (I^r)-\mu (I^{r-1}))x^r}{1-x}.$$
Comparing both expressions of the Hilbert series it follows that
$$c_j+d_j \leq r$$
 and so $$d_f\leq r-1\,.$$
In particular, $$r_J(I)=b_e.$$

Now, for regular ideals with analytic spread one we have
$$ \begin{array}{clr}
 e(F(I))& =\mu(I^r)& =e, \\
 \text{reg}(F(I))& = r_J(I)& = b_e,\\
 \mu_{F(J)}F(I))&=\mu(I^r)+f.
 \end{array}
  $$
Observe that, in this case, the reduction number $r_J(I)$ turns
out to be independent of the chosen minimal reduction $J$, as was
already noted by S. Huckaba in \cite{Huc}. Also that
$$\mu_{F(J)}(F(I))=\dim_{F(J)/F(J)_+}(F(I)/(F(J)_+F(I)))=\sum_{n=0}^{r}
\la_{A/\fm}(I^n/ \fm I^n +JI^{n-1}) .$$

In order to make a deeper analysis of the decomposition of $F(I)$
as $F(J)$-module we can rewrite it in the form
$$ F(I) \cong \bigoplus_{i=0}^r(F(J)(-i))^{\alpha _i}
\bigoplus_{i=1}^{r-1} \bigoplus_{j=1}^{r-i}
((F(J)/a^{j}F(J))(-i))^{\alpha _{i,j}} \quad \hfill( \ast \ast
).$$

Note that $\alpha_0 = 1$ and $\alpha_r\neq 0$ since $r=$ the
biggest possible degree among the generators of $F(I)$ as a graded
$F(J)$-module. Also that $$f=\displaystyle{\sum_{ \begin{array}{c}1\leq i \leq r-1\\
1\leq j \leq r-i \end{array}}} \alpha_{i,j}\,.$$

From now on we shall denote by $T(F(I))$ the $F(J)$-torsion
submodule of $F(I)$ and assume that $I$ is a regular ideal.
Observe that $T(F(I))=0$ if $r(I)=0,1$ and so in both cases $F(I)$
is a Cohen-Macaulay ring.

\begin{lem} \label{B1} Let $k,l$ and
$n$ be natural numbers. Then:
\begin{enumerate}
\item $(\fm I^{k+l}:a^l)=\fm I^k$, for $k=0$ or $k\geq r$. \item $
(\fm I^{k+1}:a)\subseteq \cdots \subseteq (\fm I^{r}:a^{r-k})=(\fm
I^{r+n}:a^{r-k+n})$, for $k=1, \ldots, r-1$.
\end{enumerate}
\end{lem}
\begin{proof}
For $k=0$ we have $(\fm I^l:a^l)=\fm $ since $a$ is analytically
independent on $I$. Now, let $k\geq r$ and $x$ be an element such
that $xa^l \in \fm I^{k+l}=a^l\fm I^k$ then $x\in \fm I^k$ since
$a^l$ is a non zero divisor in $A$.

The only non trivial inclusion in $(2)$ is $(\fm
I^{r+n}:a^{r-k+n}) \subseteq (\fm I^{r}:a^{r-k})$. Let $x$ be an
element such that $xa^{r-k+n} \in \fm I^{r+n}=a^{n}\fm I^{r}$.
Then $xa^{r-k}a^n \in a^n\fm I^r$ and now, the regularity of $a^n$
gives that $xa^{r-k} \in \fm I^r$.
\end{proof}

\begin{prop}
\label{B2} We have
$$T(F(I))=H^0_{F(I)_+}(F(I))=(0:_{F(I)} (a^0)^{r-1})=\displaystyle{\bigoplus_{k=
1}^{r-1} (I^k\cap (\fm I^{r}:a^{r-k}))}/\fm I^k .$$
\end{prop}
\begin{proof} We have that $$T(F(I))=H^0_{F(J)_+}(F(I))=H^0_{F(I)_+}(F(I))\,.$$
On the other hand,  $H^0_{F(J)_+}(F(I))= \bigcup_{l\geq 0}
(0:_{F(I)} (a^0)^l)=\bigcup_{l\geq r-1} (0:_{F(I)} (a^0)^l)$.
Thus, by the above lemma we get
$$\begin{array}{ll}
(0:_{F(I)}(a^0)^l)&=\bigoplus_{k\geq 0} (I^k\cap (\fm
I^{l+k}:a^{l}))/\fm I^k \\ &= \bigoplus_{k= 1}^{r-1} (I^k\cap (\fm
I^{l+k}:a^{l}))/\fm I^k \\ &= \bigoplus_{k= 1}^{r-1} (I^k\cap (\fm
I^{r+(l+k-r)}:a^{r-k+(l+k-r)}))/\fm I^k  \\ &= \bigoplus_{k=
1}^{r-1} (I^k\cap (\fm I^r:a^{r-k}))/\fm I^k \end{array}$$ for all
$l\geq r-1$. In particular,
$$H^0_{F(J)_+}(F(I))=(0:_{F(I)} (a^0)^{r-1})=\bigoplus_{k= 1}^{r-1}
(I^k\cap (\fm I^r:a^{r-k}))/\fm I^k \,.$$
\end{proof}

Given the two decompositions of the torsion of $F(I)$
$$T(F(I))= \bigoplus_{k=1}^{r-1} (I^k \cap (\fm I^r :a^{r-k}))/\fm I^k
\cong \bigoplus_{i=1}^{r-1} \bigoplus_{j=1}^{r-i}
((F(J)/a^jF(J))(-i))^{\alpha _{i,j}},$$ we will consider the
numbers $$ f_{k,l}:=\la ( (I^k\cap (\fm I^{k+l}:a^l))/\fm I^k).$$
Then, it is clear that $f_{k,1} \leq \cdots \leq f_{k,r-k}=\la (
[T(F(I))]_k)$ and $\la (T(F(I)))=\sum_{k=1}^{r-1} f_{k,r-k}$.
Also, that the extremal numbers $f_{k,r-k}$ are independent of the
chosen minimal reduction $J$. In addition, note that being $a$ a
non zero divisor in $A$ one has isomorphisms
$$ (I^k\cap (\fm I^{k+l}:a^l))/\fm I^k \cong (a^lI^k \cap
\fm I^{k+l})/a^l\fm I^k.$$

If $Y$ is an homogeneous element of $F(I)$ of degree $n$ we will
denote by $y$ an element of $A$ such that $Y=y^0 \in I^n/\fm I^n
\hookrightarrow F(I)$.

\begin{prop}
\label{B3} For $1\leq k \leq r-1$ and $1\leq l \leq r-k$  we have
$$f_{k,l}=\sum_{(i,j) \in \Lambda} \alpha_{i,j},$$ where $\Lambda :=\{ (i,j)
\mid 1\leq i \leq k,\, k-i+1\leq j \leq k-i+l \}$
\end{prop}
\begin{proof}
Let $\{ Y_1^{i,j},\dots ,Y_{\alpha _{i,j}}^{i,j} \}_{\{1\leq i\leq
r-1, \, 1\leq j \leq r-i\}}$ be a minimal system of homogeneous
generators of $T(F(I))$. That is,
$$T(F(I))=\bigoplus_{i=1}^{r-1}\bigoplus_{j=1}^{r-i}(F(J)Y_1^{i,j}
\oplus \cdots \oplus F(J)Y_{\alpha_{i,j}}^{i,j}) $$ with
$F(J)Y_{\ast}^{i,j} \cong (F(J)/a^jF(J))(-i)$. So,
 $$[ F(J)Y_{\ast}^{i,j} ]_k =\left\{
 \begin{array}{ll} ((a^{k-i}y_{\ast}^{i,j})+\fm I^k)/\fm I^k
 \cong A/\fm
 &\mbox{ for } i\leq k \leq i+j-1 \\ 0 &\mbox{ otherwise } \end{array} \right..$$

 Now, fixed $k$, we have $[ F(J)Y_{\ast}^{i,j} ]_k\neq 0$ for $i\leq k$
 and $j\geq k-i+1$. Moreover, $a^la^{k-i}y^{i,j}_{\ast}\in \fm I^{k+l} $ if, and only
if, $l+k-i\geq j$. Therefore, we may conclude $f_{k,l}:=\la (
(I^k\cap (\fm I^{k+l}:a^l))/\fm I^k)=\sum_{(i,j) \in \Lambda}
\alpha_{i,j}$.
\end{proof}

\begin{cor}
\label{B4} $$f= \displaystyle{\sum_{ \begin{array}{c}1\leq i \leq r-1\\
1\leq j \leq r-i \end{array}}} \alpha_{i,j} =\sum_{k=1}^{r-1}
f_{k,1}=\sum_{k=1}^{r-1} \la((aI^k\cap \fm I^{k+1})/\fm I^k),$$
and $\la (F(I)/aF(I))=\mu (I^r)+\sum_{k=1}^{r-1} \la((aI^k\cap \fm
I^{k+1})/\fm I^k)$.
\end{cor}

\begin{rem}
\label{B5} The invariants $\alpha_{i,j}$ are univocally related by
the $f_{k,l}$ (and viceversa): In fact, if we write
$$\alpha=(\alpha_{1,1},\dots,\alpha_{1,r-1},\alpha_{2,1},\dots,\alpha_{r-1,1})$$
$$F=(f_{1,1},\dots,f_{1,r-1},f_{2,1},\dots,f_{r-1,1}),$$ it is
then easy to see that there exists by proposition \ref{B3} an
invertible inferior triangular matrix $B\in M_{r(r-1)/2}$ such
that $F=B\alpha$.
\end{rem}

\begin{rem}
\label{B6} Observe that $f_{k,l}=0$ if $k\notin \{1,\dots, r-1\}$.
\end{rem}

We consider now the free part of $F(I)$ as $F(J)$-module:
$$ \begin{array}{ll}
F(I)/T(F(I))&=  A/\fm \, \displaystyle{ \bigoplus_{i=1}^{r-1}
I^i/(I^i\cap (\fm I^r:a^{r-i}))
\bigoplus_{n\geq r} I^n /\fm I^n }\\
& \cong \displaystyle{\bigoplus_{i=0}^r F(J)(-i)^{\alpha_i}}.
\end{array}$$
By convention, we put $\mu (I^0)=1$ and $\mu (I^{n})=0$ if $n<0$
for the rest of the paper.

\begin{prop}
\label{B7} For $1\leq i \leq r$ we have
$$\alpha_i=\mu(I^i)-\mu(I^{i-1})-(f_{i,r-i}-f_{i-1,r-(i-1)}).$$
\end{prop}

\begin{proof}

Put $$\alpha'_i:=\la(I^i/(I^i\cap(\fm
I^r:a^{r-i})+aI^{i-1})),\quad \alpha'_r:=\la (I^r/(\fm
I^r+aI^{r-1})),$$
 for $1\leq i\leq r-1$, and let $$\Omega=\{1,\{y_{i,1}^0,\dots ,y_{i,\alpha'_i}^0
\}\}$$ with $y_{i,j}^0 \in [F(I)/T(F(I))]_i$, for $1\leq i\leq r$
and $1\leq j \leq \alpha'_i$, such that
 $$\overline{y_{i,1}^0},\dots
,\overline{y_{i,\alpha'_i}^0} \in I^i/(I^i\cap (\fm
I^r:a^{r-i})+aI^{i-1}),\quad \overline{y_{r,1}^0},\dots
,\overline{y_{r,\alpha'_r}^0} \in I^r/(\fm I^r+aI^{r-1})
$$
form an $A/\fm $-basis. Then, $\Omega$ is a system of homogeneous
generators of $F(I)/T(F(I))$ as $F(J)$-module.

On the other hand, for $1\leq i\leq r-1$ we have the exact
sequences
$$\begin{array}{rl} 0\rightarrow aI^{i-1}/(aI^{i-1}\cap(\fm I^r:a^{r-i}))
\rightarrow & I^i/(I^i\cap(\fm I^r:a^{r-i})) \rightarrow \\
\rightarrow &I^i/(I^i\cap(\fm I^r:a^{r-i})+aI^{i-1})\rightarrow 0,
\end{array}$$ and
$$0\rightarrow (I^i\cap(\fm I^r:a^{r-i}))/\fm I^i
\rightarrow I^i/\fm I^i \rightarrow I^i/(I^i\cap(\fm
I^r:a^{r-i}))\rightarrow 0,
$$
and isomorphisms
$$aI^{i-1}/(aI^{i-1}\cap(\fm I^r:a^{r-i})) \cong
I^{i-1}/(I^{i-1}\cap(\fm I^r:a^{r-i+1})).$$ From them we obtain
$\alpha'_i=\mu(I^{i})-f_{i,r-i}-\mu(I^{i-1})+f_{i-1,r-(i-1)}$, for
$2\leq i \leq r-1$, and  $\alpha'_1=\mu(I)-f_{1,r-1}-1$.

Also we have the exact sequence
$$ 0\rightarrow aI^{r-1}/(aI^{r-1} \cap \fm I^r)\rightarrow
I^{r}/\fm I^{r} \rightarrow I^{r}/(\fm I^{r}+aI^{r-1})\rightarrow
0$$ which gives $\alpha'_r=\mu(I^{r})-\mu(I^{r-1})+f_{r-1,1}$.

Now, $\displaystyle{1+\sum_{i=1}^{r}\alpha'_i=\mu (I^r)=
\la(F(I)/T(F(I)))}$ and hence $\Omega$ is a basis of the free
$F(J)$-module $F(I)/T(F(I))$. As a consequence,
$\alpha'_i=\alpha_i$ for all $1\leq i \leq r$.
 \end{proof}

Since $\alpha_{r} \neq 0$, the following corollary extends the
invariance of the reduction number with respect to the chosen
minimal reduction.

\begin{cor}
\label{B8} The invariants $\alpha_{i}$, for $0\leq i \leq r$, are
independent of the choice of the minimal reduction.
\end{cor}
\begin{proof}
The assertion follows from Proposition \ref{B7} and Proposition
\ref{B2}.
\end{proof}

We finish this section with a lemma expressing the difference
between the minimal number of generators of the powers of $I$ in
terms of certain lengths involving minimal reductions.

\begin{lem}
\label{B9} For all $n$ and $0\leq i \leq n-1$ we have
$$\la(I^n/(\fm I^n +a^{n-i}I^i))=\mu (I^n)-\mu (I^i)+\la
((a^{n-i}I^i \cap \fm I^ n)/a^{n-i}\fm I^i).$$ In particular,
$$\la(I^n/(\fm I^n +aI^{n-1}))=\mu (I^n)-\mu (I^{n-1})+\la
((aI^{n-1} \cap \fm I^ n)/a\fm I^{n-1}).$$

\end{lem}
\begin{proof}
The exact sequence
$$0\rightarrow (a^n)/(\fm I^n \cap (a^n)) \rightarrow I^n/\fm
I^n \rightarrow I^n/(\fm I^n +(a^n)) \rightarrow 0,$$ and the
equality $\fm I^n \cap (a^n)=a^n\fm$ gives $\la (I^n/(\fm I^n
+(a^n))=\mu (I^n)-1$.

For $1\leq i \leq n-1$ we consider the exact sequences
$$0\rightarrow a^{n-i}I^i/(a^{n-i}I^i\cap \fm I^n) \rightarrow I^n/\fm
I^n \rightarrow I^n/(\fm I^n+a^{n-i}I^i) \rightarrow 0,$$
$$0\rightarrow (a^{n-i}I^i\cap \fm I^n)/a^{n-i}\fm I^i \rightarrow a^{n-i}I^i/ a^{n-i}\fm
I^i \rightarrow a^{n-i}I^i/(a^{n-i}I^i\cap \fm I^n ) \rightarrow
0$$ and the isomorphism $$a^{n-i}I^i/ a^{n-i}\fm I^i \cong I^i/\fm
I^i.$$ Then, the result follows from the additivity of $\la(\cdot
)$.
\end{proof}

\section{Buchsbaum, Cohen-Macaulay and Gorenstein properties}
\label{C}

Let $I$ be a regular ideal of $A$ with analytic spread one and
reduction number $r$ and $J=(a)\subseteq I$ be a minimal
reduction. Consider the \textit{Hilbert-Samuel} function
$$HS(F(I),n)=\la (F(I)/a^{n+1}F(I))$$
of $aF(I)$ with respect to $F(I)$. Then, $HS(F(I),n)$ is of
polynomial type of degree one and has the form
$$e_0(aF(I),F(I))(n+1) + e_1((aF(I),F(I))$$
for $n$ big enough. We shall write
$$e(aF(I),F(I)):=e_0(aF(I),F(I)).$$

In the following remark we consider in our case several well known
characterizations of the Buchsbaum property, see for instance
St{\"u}ckrad-Vogel \cite{StVo}.

\begin{rem} \label{C1}The following conditions are equivalent:
\begin{enumerate}
\item  $F(I)$ is a Buchsbaum ring. \item
$(0:_{F(I)}F(I)_+)=(0:_{F(I)} a^0)$ for any $(a)\subseteq I$
minimal reduction of $I$. \item There exists a natural number $C$
such that any $(a)\subseteq I$ minimal reduction of $I$ satisfies
$$C=\la (F(I)/aF(I))-e(aF(I),F(I)).$$ \item $F(I)_{+}\cdot
H^0_{F(I)_+}(F(I))=0$. \item $(0:a^0)=(0:(a^0)^2)$ for any minimal
reduction $(a)$ of $I$.
\end{enumerate}
In this case, $C=\la (T(F(I)))$.
\end{rem}

\begin{lem}
\label{C2} We have

\begin{enumerate}
 \item $\la(F(I)/aF(I))=\mu (I^r)+ \sum_{n=1}^{r-1} \la((aI^{n}\cap \fm
I^{n+1})/a \fm I^{n})$.

\item $\la(F(I)/a^{n+1}F(I)=\mu (I^r)(n+1)+
\sum_{k=1}^{r-1}\la((a^{r-k}I^k\cap \fm I^r)/a^{r-k}\fm I^k)$ for
all $n\geq r$.

\item $e(aF(I),F(I))=\mu (I^r)$.

\item $\la(F(I)/aF(I))- e(aF(I),F(I))= \sum_{n=1}^{r-1}
\la((aI^{n}\cap \fm I^{n+1})/a \fm I^{n})$.
\end{enumerate}
\end{lem}

\begin{proof}

(1) is proved in Corollary \ref{B4}. On the other hand, for $s\geq
1$ we have
$$\begin{array}{ll} F(I)/a^{r+s}F(I)= &
A/\fm \oplus \cdots I^r/\fm I^r \oplus  I^{r+1}/\fm I^{r+1} \oplus
\cdots \oplus I^{r+s-1}/\fm I^{r+s-1} \oplus \\ & \oplus
I^{r+s}/((a^{r+s})+\fm I^{r+s} ) \oplus I^{r+s+1}/(a^{r+s}I+\fm
I^{r+s+1} ) \oplus \cdots \\ & \cdots \oplus
I^{2r+s-1}/(a^{r+s}I^{r-1}+\fm I^{2r+s-1} ) .\end{array}$$ Observe
that there are isomorphisms
$$I^{r+i}/\fm I^{r+i} \cong I^r/\fm I^r$$
for $1\leq i \leq s-1$, and that
$$I^{r+s+i}/(a^{r+s}I^{i}+\fm
I^{r+s+i}) \cong I^r/(a^{r-i}I^i+\fm I^r)$$ for $0\leq i \leq
r-1$. Thus,
$$\la(F(I)/a^{r+s}F(I))=1+\mu (I)+\cdots +\mu (I^r)+ (s-1)\mu
(I^r)+ \sum_{i=0}^{r-1} \la(I^r/(a^{r-i}I^i+\fm I^r)).$$ Now,
$$\la(I^r/(a^{r-i}I^i+\fm I^r))=\mu (I^r)-\mu (I^i) + \la((
a^{r-i}I^i\cap \fm I^r)/a^{r-i}\fm I^i)$$ by Lemma \ref{B9} and
$$\la(F(I)/a^{r+s}F(I))=(r+s)\mu
(I^r)+\sum_{i=1}^{r-1}\la((a^{r-i}I^i\cap \fm I^r)/a^{r-i}\fm
I^i)$$ which gives (2). Finally, (3) and (4) are immediate
consequences of (1) and (2).
\end{proof}

\begin{thm}
\label{C3} Let $(A,\fm)$ be a Noeherian local ring with an
infinite residue field and let $I$ be a regular ideal of $A$ with
analytic spread one and reduction number $r$. The following
conditions are then equivalent:
\begin{enumerate}
\item $F(I)$ is a Buchsbaum ring.

\item For any minimal reduction $(a)$ of $I$ one has
$$(I^n\cap(\fm I^{n+1}:a))/\fm I^n=(I^n\cap(\fm I^{n+1}:I))/\fm
I^n$$ for $1\leq n\leq r-1$.

\item There exists an integer $C\geq 0$ such that
$$C=\sum_{n=1}^{r-1} \la ((aI^n\cap \fm I^{n+1})/a\fm I^n)$$ for any
minimal reduction $(a)$ of $I$.

\item $(0:a^0)=(0:(a^0)^{r-1})$ for any minimal reduction $(a)$ of
$I$.

\item $\la ((aI^n\cap \fm I^{ n+1})/a\fm I^n)= \la
((a^{r-n}I^n\cap \fm I^{r})/a\fm I^n)$ for any minimal reduction
$(a)$ of $I$ and $1\leq n \leq r-1$.

\item $\la((I^n\cap(\fm I^{n+1}:a))/\fm I^n)$ is independent of
the choice of the  minimal reduction $(a)$ of $I$, for $1\leq
n\leq r-1$.

\item There exists a natural number $C$ such that if $J=(a)$ is
any minimal reduction of $I$ and
$$ F(I) \cong
\bigoplus_{i=0}^r(F(J)(-i))^{\alpha _i} \bigoplus_{i=1}^{r-1}
\bigoplus_{j=1}^{r-i} ((F(J)/a^{j}F(J))(-i))^{\alpha _{i,j}}$$ is
the decomposition of $F(I)$ as $F(J)$-module, then
$$\displaystyle{C=\sum_{\begin{array}{l} 1\leq i\leq r-1 \\1\leq j
\leq r-i
\end{array} }} \alpha_{i,j}.$$

\item There exist integers $\alpha_{0} \ldots , \alpha_{r},
\alpha_{1,1} \ldots , \alpha_{r-1,1}$ such that for every
$J=(a)\subset I$ minimal reduction of $I$, the decomposition of
$F(I)$ as $F(J)$-module has the form
$$ F(I) \cong
\bigoplus_{i=0}^r(F(J)(-i))^{\alpha _i} \bigoplus_{i=1}^{r-1}
((F(J)/aF(J))(-i))^{\alpha _{i,1}}\,.$$

\end{enumerate}

\end{thm}

\begin{proof}
The equivalence $(1) \Leftrightarrow (2)$ is the corresponding one
in Remark \ref{C1}. Now, by Lemma \ref{C2} we get the equivalence
$(1) \Leftrightarrow (3)$. And by Corollary \ref{B4} we have $(3)
\Leftrightarrow (7)$.

On the other hand, $(1) \Leftrightarrow (4)$ easily follows from
$(1) \Leftrightarrow (5)$ in Remark \ref{C1}, and taking
components and their lengths in $(4)$, we get condition $(5)$.
Now, the isomorphisms $(a^{r-n}I^n\cap \fm I^{r})/a\fm I^n \cong
[H^0_{F(I)_+}(F(I))]_n$ give $(5) \Rightarrow (6)$. And it is
clear that $(6) \Rightarrow (3)$.

Finally, to get $(7) \Leftrightarrow (8)$ observe first that if we
have such a decomposition of $F(I)$ as in $(7)$, $F(I)$ is
Buchsbaum and $C=\la (T(F(I)))$ by Remark \ref{C1}. Thus
$\alpha_{i,j}=0$ for any $j\geq 2$ and, by Proposition \ref{B3},
for any $1\leq k \leq r-1$ it holds that $f_{k,
r-k}=\alpha_{k,1}$. Since the numbers $f_{k, r-k}$ are independent
of the chosen minimal reduction $J$ this implies that the
invariants $\alpha_{i,1}$, for $1\leq i \leq r-1$, are also
independent of $J$.
\end{proof}

\begin{rem}
\label{C30} Observe that as a consequence of the above theorem we have
that the Buchsbaum property of $F(I)$ is equivalent to the
invariance of the structure of $F(I)$ as a $F(J)$-module with
respect to the chosen minimal reduction $J$.
\end{rem}

Now, we recall in the following remark several characterizations
of the Cohen-Macaulay property translated to the fiber cone in
this case.

\begin{rem} \label{C4}The following conditions are equivalent:
\begin{enumerate}
\item $F(I)$ is a Cohen-Macaulay  ring.
 \item $(0:_{F(I)} a^0)=0$ for
all (some) $(a)\subseteq I$ minimal reduction. \item $\la
(F(I)/aF(I))-e(aF(I),F(I))=0$ for every (some) $(a)\subseteq I$
minimal reduction.
 \item $H^0_{F(I)_+}(F(I))=0$. \item
$F(I)$ is a free $F(J)$-module, for every (some) minimal reduction
$J$ of $I$.
\end{enumerate}
\end{rem}

\begin{thm}
\label{C5} Let $(A,\fm)$ be a Noetherian local ring with an
infinite residue field and let $I$ be a regular ideal of $A$ with
analytic spread one and reduction number $r$. The following
conditions are then equivalent:

\begin{enumerate}
\item $F(I)$ is a Cohen-Macaulay ring.

\item For every (some) minimal reduction $(a)$ of $I$,
$I^n\cap(\fm I^{n+1}:a)=\fm I^n$ for $1\leq n\leq r-1$.

\item For every (some) minimal reduction $(a)$ of $I$, $\la
(F(I)/aF(I))=\mu (I^r)$.

\item For every (some) minimal reduction $(a)$ of $I$, $I^k\cap
(\fm I^r:a^{r-k})=\fm I^k$ for $1\leq k \leq r-1$.

\item For every (some) minimal reduction $(a)$ of $I$, $aI^n\cap
\fm I^{n+1}=a\fm I^n$ for $1\leq n\leq r-1$.

 \item For every (some) minimal reduction $(a)$ of $I$, $\la (I^n/(\fm
I^n+aI^{n-1}))=\mu (I^n)-\mu (I^{n-1})$ for $1\leq n \leq r$.

\item For every (some) $J=(a)$ minimal reduction of $I$ the
decomposition of $F(I)$ as $F(J)$ module has the form
$$\begin{array}{ll} F(I) &\cong  F(J)\bigoplus (F(J)(-1))^{\mu
(I)
-1} \oplus \cdots \oplus (F(J)(-r))^{\mu (I^r) - \mu (I^{r-1})} \\
& = \bigoplus_{i=0}^{r} F(J)(-i)^{\mu (I^i)-\mu
(I^{i-1})}.\end{array}$$

\end{enumerate}
\end{thm}

\begin{proof}
Sentences $(1),(2),(3)$ and $(4)$ correspond to the same ones in
Remark \ref{C4} and so the equivalences. The equivalence
$(2)\Leftrightarrow (5)$ follows for the regularity of $a$ in $A$,
while Lemma \ref{B9} gives $(5)\Leftrightarrow (6)$.

On the other hand, If there exists a minimal reduction $J$ of $I$
as in $(7)$ then $F(I)$ is a free $F(J)$-module and so $F(I)$ is a
Cohen-Macaulay ring. Conversely, if  $F(I)$ is Cohen-Macaulay then
it is a free $F(J)$-module for every $J$ minimal reduction, and it
has a decomposition as a direct sum of simple free $F(J)$-modules
$$F(I) \cong \bigoplus_{i=0}^{r-1} F(I)(-i)^{\alpha_i} $$
where, by Proposition \ref{B7}, $\alpha_i=\mu (I^i)-\mu (I^{i-1})
-f_{i,r-i}+f_{i-1,r-i+1}$ and $f_{k,r-k}=\la((I^k\cap (\fm
I^r:a^{r-k})/\fm I^k)=0$ for $1\leq k \leq r-1$.
\end{proof}

\begin{cor}
\label{C10} Let $(A,\fm)$ be a Noeherian local ring with an
infinite residue field and let $I$ be a regular ideal of $A$ with
analytic spread one and reduction number $r$. If $F(I)$ is
Cohen-Macaulay then $\mu (I^{n-1}) <\mu (I^{n})$ for $1\leq n \leq
r$ and the postulation number $fp(I)=r-1$.
\end{cor}
\begin{proof} Let $J=(a)$ be a minimal reduction of $I$.
If $F(I)$ is Cohen-Macaulay then, for $1\leq n \leq r$, we have by
the Nakayama Lemma and the previous result that $0<\la (I^n/(\fm
I^n+aI^{n-1}))=\mu (I^n)-\mu (I^{n-1})$. Thus, $1<\mu(I) <\mu
(I^2)< \cdots  <\mu (I^{r-1}) <e(F(I))=\mu (I^r)=\mu (I^m)$ for
all $m\geq r$. (See also D'Anna-Guerrieri-Heinzer
\cite[Proposition 3.2]{D'AnGuHe1} for the computation of
postulation number.)\end{proof}

To conclude this section we study the Gorenstein property of
$F(I)$.

\begin{rem} \label{C6}The following conditions are equivalent:
\begin{enumerate}
\item $F(I)$ is a Gorenstein  ring.
 \item $F(I)$ is a Cohen-Macaulay ring and $\mathrm{type}(F(I))=1$
 \item For every (some) minimal reduction $J$ of $I$, $F(I)$ is a free $F(J)$-module
 and and the canonical module $\omega _{F(I)}\cong F(I)(k)$, for some $k\in \mathbb Z$.
\end{enumerate}
\end{rem}

\begin{lem}
\label{C7} Let $(A,\fm)$ be a Noeherian local ring with an
infinite residue field and let $I$ be a regular ideal of $A$ with
analytic spread one and reduction number $r$. Assume that $F(I)$
is a Cohen-Macaulay ring and let $J=(a)\subseteq I$ be a minimal
reduction. Then
$$ \mathrm{type}(F(I))=\sum _{i=1}^{r-1}\la( (I^i\cap (aI^i+\fm I^{i+1}:
I))/(aI^{i-1}+\fm I^i)) + \la( I^r/(aI^{r-1}+\fm I^r)).$$
\end{lem}

\begin{proof}
Let $(a)$ be a minimal reduction of $I$. Since $F(I)$ is
Cohen-Macaulay, $a^0$ is a regular element in $F(I)$ and
$$\mathrm{type}(F(I))=\la ( \text{Ext}^1(A/\fm ,F(I)))= \la
(\text{Socle}( F(I)/aF(I))) = \la ((0:_{F'}F'_+))$$ where $
F'=F(I)/aF(I)$. The statement then follows from the equality
$$(0:_{F'}F'_+)= \bigoplus _{i=1}^{r-1} ((I^i\cap (aI^i+\fm
I^{i+1}: I))/(aI^{i-1}+\fm I^i)) \oplus  I^r/(aI^{r-1}+\fm I^r).$$
\end{proof}

Assume that $F(I)$ is Cohen-Macaulay. In the following lemma we
describe the structure as a $F(J)$-module of the canonical module
of $F(I)$. Recall that, since $F(J)$ is a polynomial ring in one
variable over a field, the a-invariant of $F(J)$ is $-1$ and
$\omega_{F(J)}\simeq F(J)(-1)$.

\begin{lem}
\label{C8} Let $(A,\fm)$ be a Noetherian local ring with an
infinite residue field and let $I$ be a regular ideal of $A$ with
analytic spread one and reduction number $r$. Assume that $F(I)$
is a Cohen-Macaulay ring and let $J$ be a minimal reduction of
$I$. Then,
$$\omega_{F(I)} \simeq \text{Hom}_{F(J)}(F(I),F(J)(-1)) \cong \bigoplus_{i=0}^{r}(
F(J)(i-1))^{\mu (I^{i})-\mu (I^{i-1})}$$ and the $a$-invariant of
$F(I)$ is $r-1$.
\end{lem}

\begin{proof}
We may write $F(I)\cong \bigoplus_{i=0}^{r}( F(J)(-i))^{\mu
(I^{i})-\mu (I^{i-1})}$. Then, by local duality,
$$ \begin{array}{ll}
\omega_{F(I)} & \cong \text{Hom}_{F(J)}(F(I),F(J)(-1)) \\ & \cong
\text{Hom}_{F(J)}( \bigoplus_{i=0}^{r}( F(J)(-i))^{\mu (I^{i})-\mu
(I^{i-1})},
F(J)(-1)) \\
& \cong \bigoplus_{i=0}^{r}\text{Hom}_{F(J)}(( F(J)(-i))^{\mu
(I^{i})-\mu (I^{i-1})}, F(J)(-1))\\
& \cong \bigoplus_{i=0}^{r}(\text{Hom}_{F(J)}(F(J)(-i),F(J)(-1)
)^{\mu (I^i)- \mu(I^{i-1})} \\
& \cong
\bigoplus_{i=0}^{r}(\text{Hom}_{F(J)}(F(J),F(J))(i-1))^{\mu
(I^i)- \mu(I^{i-1})} \\
& \cong \bigoplus_{i=0}^{r}(F(J)(i-1))^{\mu (I^i)- \mu(I^{i-1})}
 .\end{array}
$$
\end{proof}

\begin{thm}
\label{C9} Let $(A,\fm)$ be a Noeherian local ring with an
infinite residue field and let $I$ be a regular ideal of $A$ with
analytic spread one and reduction number $r$. Then, the following
conditions are equivalent:
\begin{enumerate}
\item $F(I)$ is a Gorenstein ring.
 \item  $\mu (I^r)=\mu (I^{r-1})+1$, and
for every (some) minimal reduction $J=(a)$ of $I$ the following
equalities hold
$$ I^n\cap(\fm I^{n+1}:a)=\fm I^n \text{ and }
I^n\cap(aI^{n}+\fm I^{n+1}:I)=aI^{n-1}+\fm I^n $$ for $1\leq n\leq
r-1$
\end{enumerate}
In this case, the decomposition of $F(I)$ as the direct sum of
cyclic $F(J)$-module has the form
$$\begin{array}{ll}F(I) & \cong \bigoplus_{i=0}^{r} F(J)(-i)^{\mu (I^i)-\mu
(I^{i-1})},  \end{array}$$ with $\mu (I^i)-\mu (I^{i-1})=\mu
(I^{r-i})-\mu (I^{r-i-1})$ for $0\leq i \leq r$.
\end{thm}

\begin{proof}
Let $J=(a)\subseteq I$ be a minimal reduction. We know that $F(I)$
is Cohen-Macaulay if and only if $ I^n\cap(\fm I^{n+1}:a)=\fm I^n
$ for $1\leq n \leq r-1$. Moreover, in this case $\la(
I^r/(aI^{r-1}+\fm I^r))=\mu (I^r)-\mu (I^{r-1})>0$. Thus, we may
assume that  $F(I)$ is Cohen-Macaulay and so it is Gorenstein if,
and only if, $\mathrm{type}(F(I))=1$. By Lemma \ref{C7} this is
equivalent to $1=\la( I^r/(aI^{r-1}+\fm I^r))=\mu (I^r)-\mu
(I^{r-1})$ and $I^n\cap(aI^{n}+\fm I^{n+1}:I)=aI^{n-1}+\fm I^n $,
for $1\leq n\leq r-1$.

Now, by Theorem \ref{C5} $F(I)$ is Cohen-Macaulay if, and only if,
$$F(I)\simeq  \bigoplus_{i=0}^{r} F(J)(-i)^{\mu (I^i)-\mu (I^{i-1})}.$$
And by Lemma \ref{C8}, $$\omega _{F(I)} \cong
\text{Hom}_{F(J)}(F(I), F(J)(-1)) \cong \bigoplus_{i=0}^{r}(
F(J)(i-1))^{\mu (I^{i})-\mu (I^{i-1})}$$ with $a(F(I))=r-1$. On
the other hand $F(I)(k)\cong \bigoplus_{i=0}^{r}(F(J)(k-i))^{\mu
(I^{i}) -\mu (I^{i-1})} $. Hence, if $F(I)$ is Gorenstein, $\omega
_{F(I)} \cong F(I)(r-1 )$ and just comparing we get $\mu
(I^{i})-\mu (I^{i-1}) =\mu (I^{r-i})-\mu (I^{r-i-1})$. (Observe
that one may also obtain these equalities from the well known fact
that the $h$-vector of a Gorenstein graded algebra is symmetric.)
\end{proof}

\section{Applications and examples}
\label{D}

We may first apply the results in the above section to the case of
ideals with small reduction number.

\begin{prop}
\label{D1} Let $(A,\fm)$ be a Noeherian local ring with an
infinite residue field and let $I$ be a regular ideal of $A$ with
analytic spread one and reduction number $1$. Then,
\begin{enumerate}
\item For any $J$ minimal reduction of $I$, $F(I)\cong
F(J)\oplus (F(J)(-1))^{\mu (I)-1}$.
 \item $F(I)$ is  a Cohen-Macaulay ring.
 \item $F(I)$ is a Gorenstein ring if and only if $\mu (I)=2$.
 \end{enumerate}
 \end{prop}
 \begin{proof}
We have already noted in Section \ref{B} that $F(I)$ is
Cohen-Macaulay if $r(I)=1$. Then, apply Theorem \ref{B5} to get
(1). Finally, (3) is a consequence of Theorem \ref{C9}.
\end{proof}

\begin{prop}
\label{D2} Let $(A,\fm)$ be a Noeherian local ring with an
infinite residue field and let $I$ be a regular ideal of $A$ with
analytic spread one and reduction number $2$. Then,
\begin{enumerate}
\item For any minimal reduction $J=(a)$ of $I$,
$$F(I)\cong F(J)\oplus (F(J)(-1))^{\mu (I)-1 -\alpha} \oplus
(F(J)(-2))^{\mu (I^2)-\mu (I) +\alpha } \oplus
(F(J)/JF(J)(-1))^{\alpha }, $$ where $\alpha =\la ((aI \cap \fm
I^2)/a\fm I)$.
 \item $F(I)$ is  a Buchsbaum  ring.
 \item The following condtions are equivalent:
 \begin{enumerate}
 \item $F(I)$ is a Cohen-Macaulay ring.
 \item  $aI \cap \fm I^2=a\fm I$ for every (some) minimal reduction $(a)$ of $I$.
\item $\mu (I^2)-\mu (I)=\la (I^2/(\fm I^2+aI))$ for every (some)
minimal reduction $(a)$ of $I$.
\end{enumerate}
 \item
The following conditions are equivalent:
 \begin{enumerate}
 \item $F(I)$ is a Gorenstein ring.
 \item $\mu( I^2)- \mu (I)=\la (I^2/(\fm I^2+aI))=1$ for every (some)
minimal reduction $(a)$ of $I$ and $I\cap (aI+\fm I^2:I)=(a)+\fm
I$.

 \end{enumerate}
\end{enumerate}
 \end{prop}
 \begin{proof}
Let $J=(a)$ be a minimal reduction of $I$. Then, (1) is a a direct
consequence of Corollary \ref{B4} and Proposition \ref{B7}. By
Theorem \ref{C3}, $F(I)$ is Buchsbaum if and only if $I\cap (\fm
I^2:a)=I\cap (\fm I^2:I)$ for every minimal reduction $(a)$ of
$I$. Let $x\in I$ such that $xa \in \fm I^2$. Then, for any $y\in
I$ one has $axy \in \fm I^3=a\fm I^2$ and so $xy\in \fm I^2$ since
$a$ is regular, that is, $x \in I\cap (\fm I^2:I)$. In fact, note
that $\alpha=\la (T(F(I))) $ and so, by Proposition \ref{B2}, it
is independent of the choice of $J$, which also proves the
Buchsbaum property of $F(I)$ by Theorem \ref{C3}. Now, (3) and (4)
a are direct consequence of Theorem \ref{C5} and Theorem \ref{C9},
respectively, for $r=2$.
\end{proof}

Assume in addition that $I=\fm$. Then, condition (3), (b) in the
above proposition trivially holds since $\fm ^3 = a\fm^2$. One can
also see in this case that condition (4), (b) holds if $A$ is
Gorenstein (see, for instance, Proposition 3.3 and the final part
of the proof of Theorem 3.4 in J. Sally \cite{Sa}). As for the
Buchsbaum property, it is known that $F(I)=G(\fm )$ is Buchsbaum
if $A$ is Buchsbaum of $\dim A=1$ and $r(\fm) =2$, see S. Goto
\cite[Proposition 7.4]{Go}.

It is easy to see that if $I=\fm$ and $r(\fm)=3$, then condition
(2) in Theorem \ref{C3} holds and so $F(I)=G(\fm )$ is Buchsbaum
in this case, see also S. Goto \cite[Proposition 7.7]{Go}.
Nevertheless, this result cannot be extended to more general fiber
cones as the following examples will show:

D'Anna-Guerrieri-Heinzer describe in \cite[Example 2.3]{D'AnGuHe1}
a family $(R_n,\fm _n)$ (for $n\geq 3$) of one-dimensional local
Cohen-Macaulay rings and $\fm _n$-primary ideals $I_n$ for which
$\mu (I_n)=n$, $r(I_n)=n-1$ and the fiber cone $F(I_n)$ is not
Cohen-Macaulay. Moreover, $\mu ((I_n)^j)=\mu (I_n)$, for all
$j\geq 1$. For our purposes, we are going to consider the
particular cases $n=3,4$.

\begin{exam}
 Consider $R_3=K[[t^6,t^{11},t^{15},t^{31}]]$ and
$I_3=I=(t^6,t^{11},t^{31})$. Since $r(I)=2$, the fiber cone $F(I)$
is a Buchsbaum ring (and not Cohen-Macaulay). For any minimal
reduction $J=(a)$ of $I$, the structure of $F(I)$ as $F(J)$-module
is, by Proposition \ref{D2},
$$F(I)\cong F(J) \oplus F(J)(-1)\oplus F(J)(-2) \oplus
(F(J)/JF(J))(-1),$$ since in this case $\mu (I)=\mu (I^2)=3$, and
$\alpha =\la ((aI\cap \fm I^2)/\la \fm I)=\la (I^2/(\fm
I^2+aI))=\la ( (t^{12},t^{17},t^{22})/\fm I^2 +(t^{12},t^{17},
t^{37}))=1$.
\end{exam}

\begin{exam}
Consider $R_4=K[[t^8,t^{15},t^{28},t^{50},t^{57}]]$ and
$I_4=I=(t^8,t^{15},t^{50},t^{57})$. We claim that $F(I)$ is not a
Buchsbaum ring.

By \cite{D'AnGuHe1}), $(t^8)$ is a minimal reduction of $I$, $\mu
(I)=4$ and $r(I)=3$. In order to prove that $F(I)$ is not
Buchsbaum we will show that there exists an element $x\in
I\setminus \fm I$ such that $x\cdot (t^8)^2 \in \fm I^3$ and
$x\cdot t^8 \notin \fm I^2$; This implies that $x^0\in
(0:({t^8}^0)^2)$, $0\neq {t^8}^0\cdot x^0 =\overline{t^8x} \in
I^2/\fm I^2 \subset F(I)$ and so $F(I)_+ \cdot (0:({t^8}^0)^2)\neq
0$. Take $x=t^{57}$: Then $x$ fulfills the conditions since
$t^{16} t^{57}=t^{73}=t^{28}(t^{15})^3 \in \fm I^3$ and $t^8
t^{57}=t^{65} \notin \fm I^2$.

In order to describe the structure of $F(I)$ as $F(t^8)$-module we
observe that, since $\mu (I^4)=\mu (I)= 4 $ and $f_{11}<f_{12}$,
$$
\begin{array}{l} f_{11}=\la (I^2/(\fm
I^2+t^8I)=\la((t^{16},t^{23},t^{30},t^{65})/(\fm I^2
+(t^{16},t^{23},t^{58},t^{65}))=1,\\
f_{12}=\la (I^3/(\fm I^3+t^{16}I)=\la
((t^{24},t^{31},t^{45},t^{38})/\fm
I^3+(t^{24},t^{31},t^{66},t^{73}))=2,\\
 f_{21}=\la (I^3/(\fm I^3+t^8I^2)=\la((t^{24},t^{31},t^{45},t^{38})/(\fm
 I^3
+(t^{24},t^{31},t^{38},t^{73}))=1. \end{array}$$ So,
$\alpha_1=\mu(I)-1-f_{12}=1$, $\alpha_2=f_{12}-f_{21}=1$,
$\alpha_3=f_{21}=1$,
$\alpha_{11}=f_{11}=1$,$\alpha_{12}=f_{12}-f_{11}=1$,
$\alpha_{21}=f_{21}-f_{12}+f_{11}=0$ and
$$\begin{array}{ll
  } F(I)\cong & F(t^8)\oplus F(t^8)(-1)\oplus F(t^8)(-2)\oplus
F(t^8)(-3)\oplus \\ & \oplus (F(t^8)/(t^8F(t^8))(-1)) \oplus
(F(t^8)/(t^{16}F(t^8))(-1)).\end{array}$$

On the other hand, $t^{57}I^3 \subseteq \fm I^4$. Thus, writing
$t^8x=(t^8+t^{57})x-t^{57}x$ for any $x\in I^3$ one gets that
$I^4=t^8I^3 \subseteq (t^8+t^{57})I^3 +\fm I^4$. By Nakayama's
Lemma we obtain that $(t^8+t^{57})$ is also a minimal reduction of
$I$ such that $$\la ((I\cap(\fm I^2:(t^8+t^{57})))/\fm I)=\la
(I^2/(\fm I^2 +(t^6+t^{57})I)=2.$$ In this case, the structure of
of $F(I)$ as $F(t^8+t^{57})$-module is given by
$$\begin{array}{ll
  } F(I)\cong & F(t^8+t^{57})\oplus F(t^8+t^{57})(-1)\oplus F(t^8+t^{57})(-2)\oplus
F(t^8+t^{57})(-3)\oplus \\ & \oplus ((F(t^8)/(t^8F(t^8))(-1)))^2
\oplus (F(t^8)/(t^{8}F(t^8))(-2)).\end{array}$$
\end{exam}

This shows that the structure of $F(I)$ as $F(J)$-module may
depend on the chosen reduction $J$ when $F(I)$ is not Buchsbaum.

In the next lemma we prove a closed formula for the minimal number
of generators of the powers $I^n$ of a regular ideal with analytic
spread one. It also provides an easy proof in this case of a well
known of Eakin and Sathaye \cite{EaSa}, see \cite[9.39]{Va} for a
general proof, Hoa-Trung \cite{HoTr} for a combinatorial approach,
or the more recent proof by G. Caviglia \cite{Ca}.

\begin{lem}
\label{D3} Let $(A,\fm)$ be a Noeherian local ring with an
infinite residue field and let $I$ be a regular ideal of $A$ with
analytic spread one and reduction number $r$. Let $(a)$ be a
minimal reduction if $I$. Then
$$\mu(I^{n})=1+\sum_{i=2}^{n} \la \left ( \frac{\fm
I^{n}+a^{i-1}I^{n-i+1}}{\fm I^{n}+a^{i}I^{n-i}}\right ) +\la \left
( \frac{I^{n}}{\fm I^{n}+aI^{n-1}} \right ),$$ and, for $1\leq n
\leq r$, we have $\mu (I^{n})\geq n+1$.
\end{lem}

\begin{proof}

Put $J=(a)$. Fixed $n$, we consider the following exact sequences
$$
\begin{array}{c}
0\rightarrow J^{n}/(\fm I^{n}\cap J^{n}) \rightarrow I^{n}/\fm
I^{n} \rightarrow I^{n}/(\fm I^{n}+J^{n}) \rightarrow 0 ,\\

0\rightarrow (\fm I^{n}+J^{n-1}I)/(\fm I^{n}+ J^{n}) \rightarrow
I^{n}/(\fm
I^{n}+J^{n}) \rightarrow I^{n}/(\fm I^{n}+J^{n-1}I) \rightarrow 0 ,\\

0\rightarrow (\fm I^{n}+J^{n-2}I^2)/(\fm I^{n}+ J^{n-1}I)
\rightarrow I^{n}/(\fm
I^{n}+J^{n-1}I) \rightarrow I^{n}/(\fm I^{n}+J^{n-2}I^2) \rightarrow 0 ,\\

\cdots \cdots \\

0\rightarrow (\fm I^{n}+JI^{n-1})/(\fm I^{n}+ J^{2}I^{n-2})
\rightarrow I^{n}/(\fm I^{n}+J^{2}I^{n-2}) \rightarrow I^{n}/(\fm
I^{n}+JI^{n-1}) \rightarrow 0 .
\end{array}$$

Then
$$\mu(I^{n})=1+\sum_{i=2}^{n} \la \left ( \frac{\fm
I^{n}+J^{i-1}I^{n-i+1}}{\fm I^{n}+J^{i}I^{n-i}}\right ) +\la \left
( \frac{I^{n}}{\fm I^{n}+JI^{n-1}}\right ).$$

Let now $n\leq r$. It is clear that $\la (I^{n}/(\fm
I^{n}+JI^{n-1}))>0$ by the Nakayama's lemma.

Claim: $\displaystyle{\la \left (\frac{\fm
I^{n}+J^{i-1}I^{n-i+1}}{\fm I^{n}+J^{i}I^{n-i}}\right )>0.}$

Assume that there exist $n$ and $i$, $2\leq i \leq n \leq r$, such
that $\fm I^{n}+J^{i-1}I^{n-i+1} \subseteq \fm I^{n}
+J^{i}I^{n-i}$. Then $J^{i-1}I^{n-i+1} \subseteq \fm I^{n}
+J^{i}I^{n-i}$, and for all $k\geq n-i+1$ we also have
(multiplying by $I^{k-(n-i+1)})$ that $J^{i-1}I^{k} \subseteq \fm
I^{k+i-1}+J^{i}I^{k-1}$. In particular, since $r\geq n-i+1$,
$I^{r+i-1}=J^{i-1}I^{r} \subseteq \fm I^{r+i-1}+J^{i}I^{r-1}$. The
Nakayama's lemma implies now $I^{r+i-1}=J^{i}I^{r-1}$. On the
other hand $I^{r+i-1}=J^{i-1}I^{r}$ $(i-1\geq 1)$. Therefore
$a^{i-1}I^{r}=a^{i}I^{r-1}$ with $a$  regular, and so
$I^{r}=aI^{r-1}$ which contradicts the definition of $r$.
\end{proof}

Now we shall apply the above lemma to the case of ideals generated
by exactly two elements.

\begin{prop}
\label{D4} Let $(A,\fm)$ be a Noetherian local ring with an
infinite residue field and let $I$ be a regular ideal with
analytic spread one. Assume that $I$ is minimally generated by $2$
elements. Let $J=(a)$ be a minimal reduction of $I$. Then, $F(I)
\cong \displaystyle{\bigoplus_{i=0}^rF(J)(-i)}$ and $F(I)$ is a
Gorenstein ring.
\end{prop}

\begin{proof}
Let $(a)$ be a minimal reduction of $I$. Since $\mu (I)=2$, it's
easy to see that $\mu (I^n)\leq n+1$ for any $n\geq 1$. Thus, by
the above lemma we have
$$\mu(I^{n})=1+\sum_{i=2}^{n} \la \left ( \frac{\fm
I^{n}+a^{i-1}I^{n-i+1}}{\fm I^{n}+a^{i}I^{n-i}}\right ) +\la \left
( \frac{I^{n}}{\fm I^{n}+aI^{n-1}} \right )= n+1$$ for $1\leq
n\leq r$.  This implies
$$\la \left ( \frac{\fm I^{n}+a^{i-1}I^{n-i+1}}{\fm
I^{n}+a^{i}I^{n-i}}\right ) =\la \left ( \frac{I^{n}}{\fm
I^{n}+aI^{n-1}}\right )=1$$ for any $1\leq n\leq r$ and $2\leq
i\leq n$.

Furthermore, by Lemma \ref{B9}
$$ \la ( I^{n}/\fm
I^{n}+aI^{n-1})=\mu (I^n)-\mu (I^{n-1}) +\la ((aI^{n-1}\cap \fm
I^{n}) /a\fm I^{n-1}). $$ Thus $\la ((aI^{n-1}\cap \fm I^{n})
/a\fm I^{n-1})=0$, and
$$F(I)\cong F(J) \oplus F(J)(-1) \cdots  \oplus F(J)(-r)
=\bigoplus_{i=0}^rF(J)(-i).$$ On the other hand, $$\la (I^n\cap
(aI^n+\fm I^{n+1}:I)/aI^{n-1}+\fm I^n)\leq \la (I^n/aI^{n-1}+\fm
I^n)=1$$ for any $1\leq n\leq r-1$. Thus $(aI^n+\fm I^{n+1}:I)\neq
aI^{n-1}+\fm I^n$ if, and only if, $\la (I^n\cap (aI^n+\fm
I^{n+1}:I)/aI^{n-1}+\fm I^n)=1$ if, and only if, $I^n\cap
(aI^n+\fm I^{n+1}:I)=I^n$ if, and only if, $I^n\subset (aI^n+\fm
I^{n+1}:I)$ if, and only if, $I^{n+1}\subset aI^n+\fm I^{n+1}$ if,
and only if, $I^{n+1}\subset aI^n$, which is not possible for
$n\leq r-1$. Thus, by Theorem \ref{C9}, $F(I)$ is Gorenstein.
\end{proof}

In fact, it is proven in D'Anna-Guerrieri-Heinzer
\cite[Proposition 3.5]{D'AnGuHe1} that the fiber cone of a regular
ideal minimally generated by two elements having a principal
reduction is a complete intersection. This result has been
extended by Heinzer-Kim \cite[Theorem 5.6]{HeKi} to ideals of
arbitrary analytic analytic spread $l>0$, minimally generated by
$l+1$ elements and having a minimal reduction generated by a
regular sequence, such that the associated graded ring has a
homogeneous regular sequence of length at least $l-1$, see also
Jayanthan-Puthenpurakal-Verma \cite[Proposition 4.2]{JaPuVe}.

Assume now that $A$ is Cohen-Macaulay of dimension 1 and let $I$
be an $\fm$-primary ideal. Then
$$\la(I^{n+1}/aI^n)=\la(A/(a))-\la(I^n/I^{n+1})$$ for all
$(a)\subseteq I$, minimal reduction of $I$. Since
$\la(A/(a))=e(I)$, the multiplicity of the ideal $I$, the lengths
$\la(I^{n+1}/aI^n)$ are independent of $(a)$. (An $\fm$-primary
ideal $I$ in a Cohen-Macaulay ring $A$ such that $\la (I^2/JI)=1$
for any minimal reduction $J$ of $I$ is called a Sally ideal in
\cite{JaPuVe}.)

\begin{prop}
\label{D5} (See also \cite[Theorem 3.3]{JaPuVe}) Let $(A,\fm)$ be
a Cohen-Macaulay local ring of dimension $1$ with an infinite
residue field and let $I$ be an $\fm$-primary ideal with reduction
number $r$. Let $J=(a)\subset I$ be a minimal reduction of $I$ and
Assume that $e(I)-\la(I/I^2)=1$. Then
\begin{enumerate}
\item $F(I)$ is a Buchsbaum ring and $F(I)\simeq$
$$F(J)\oplus F(J)(-1)^{\mu (I^2)-2}\bigoplus_{i=2}^{r-1}
F(J)(-i)^{\mu (I^{i+1})-\mu (I^i)} \oplus F(J)(-r) $$
$$\bigoplus_{i=1}^{r-1}((F(J)/aF(J))(-i))^{\mu(I^i)-\mu(I^{i+1})+1}\,.$$
\item The following conditions are equivalent
\begin{enumerate}
\item $F(I)$ is a Cohen-Macaulay ring.

\item  $\mu (I^{n+1})=\mu(I^n)+1$ for all $1\leq n \leq r-1$

\item $\mu (I^2)=\mu (I)+1$.

\item $\fm I^2=a\fm I$ for every (some) $(a)$ minimal reduction of
$I$.

\item For every (some) $J$ minimal reduction of $I$ there exists
an isomorphism $$F(I)\cong F(J)\oplus F(J)(-1)^{\mu (I)-1} \oplus
F(J)(-2) \oplus \cdots \oplus F(J)(-r).$$
\end{enumerate}

In this case, $\text{type} (F(I))=\la ((I\cap (aI :I))/((a)+\fm
I))+1$.

\item If $r\geq 3$, $F(I)$ is Gorenstein if and only if $\mu
(I)=2$.

\end{enumerate}
\end{prop}

\begin{proof}
Let $(a)$ be a minimal reduction of $I$. Then $\la (I^2/aI)=1$.
This condition implies that $\fm I^{n+1} \subseteq aI^n$ for all
$n$ and $\la (I^{n+1}/aI^n)=1$ for all $1\leq n <r$ and so,
applying Lemma \ref{B9} we obtain for $1\leq n \leq r-1$ the
equalities $$\la (\fm I^{n+1}/a \fm I^n)= \la ((aI^n\cap \fm
I^{n+1})/a\fm I^n)= 1+\mu (I^n)-\mu (I^{n+1}).
$$
Therefore, $\la ((aI^n\cap \fm I^{n+1})/a\fm I^n)$ is independent
of the reduction and by Theorem \ref{C3} the ring $F(I)$ is
Buchsbaum. Moreover, $\alpha_{i,j}=0$ for $j\geq 2$ and by
Proposition \ref{B3} $\alpha_{k,1}=f_{k,r-k}= \cdots = f_{k,1}=
\la ((aI^n\cap \fm I^{n+1})/a\fm I^n)$ for any $1\leq k\leq r-1$.
Then, by Proposition \ref{B7} we may get the the values of
$\alpha_i$'s for $0\leq i \leq r$.

On the other hand, by Theorem \ref{C5} $F(I)$ is Cohen-Macaulay if
and only if for $1\leq n \leq r-1$ one has $\la (\fm I^{n+1}/a \fm
I^n)=0$ and (2) follows easily.

Assume now that $r\geq 3$ and $F(I)$ is Gorenstein. Then, by
Theorem \ref{C9} $\mu (I)-1=\mu (I^{r-1})-\mu (I^{r-2})=1$ and so
$\mu(I)=2$. Finally, by Proposition \ref{D4} we have the converse.
\end{proof}

In \cite[6]{JaPuVe} one may find various interesting examples of
Sally ideals with reduction number two. The following one is a
Sally ideal of reduction number three whose fiber cone is not
Cohen-Macaulay, see also \cite[Example 2.2]{Sa}.

Let $A=k[[t^4, t^5, t^{11}]]$ where $k$ is any field and $t$ an
indeterminate. Let $I= \fm =(t^4, t^5, t^{11})$ be the maximal
ideal of $A$. Then, one can easily see that $J=(t^4)$ is a minimal
reduction of $I$, the reduction number of $I$ is $3$, and $\la
(I^2/t^4I)=1$. Moreover, $\mu (I)=\mu (I^2)=3$ and $\mu(I^3)=4$.
Hence $$F(I)\simeq F(J)\oplus F(J)(-1)\oplus F(J)(-2) \oplus
F(J)(-3) \oplus F(J)/t^4F(J)(-1).$$

\section{Ideals of higher analytic spread }
\label{E}

In this section we give some applications to ideals of higher
analytic spread.

Let $I$ be an ideal of $A$. We will denote by $G(I)$ the
associated graded ring of $I$. Given $a\in I$ we will set
$a^{\ast} \in I/I^2 \hookrightarrow G(I)$ and $a^0 \in I/\fm
I\hookrightarrow F(I)$. Let $a_1, \ldots, a_k$ be a family of
elements in $I$. Then, by the well known Valabrega-Valla criterium
$a_1^{\ast}, \ldots, a_k^{\ast}$ is a regular sequence in $G(I)$
if and only if
\begin{enumerate}
\item $a_1,\dots , a_k$ is a regular sequence in $A$. \item
$(a_1,\dots ,a_{k})\cap I^n =(a_1,\dots ,a_{k})I^{n-1}$, for all
$n\geq 0$.
\end{enumerate}
In this case, there are natural isomorphisms
$$G(I/(a_1,\ldots ,a_k))\simeq G(I)/(a_1^{\ast}, \ldots
,a_k^{\ast})$$ $$F(I/(a_1, \ldots ,a_k))\simeq F(I)/(a_1^0, \ldots
,a_k^0)$$

Let $S$ be a standard $\mathbb N$-graded algebra over a local
ring. Recall that a sequence of homogeneous elements $x_1,\dots
,x_k$ is called {\it filter-regular} if for any $1\leq i\leq k$,
$[(x_1,\dots,x_{i-1}):x_i]_n =[(x_1,\dots ,x_{i-1})]_n$ for $n\gg
0$.

Let $a\in I$. If $a\in I\setminus I^2$, $a$ is a superficial
element for $I$ if and only if $a^*$ is filter regular in $G(I)$
(see for instance \cite[Lemma 2.3]{Co}). In analogy to this
situation, Jayanthan-Verma \cite{JaVe2} define $a^0\neq 0$ to be
superficial in $F(I)$ if and only if $a^0$ is filter-regular in
$F(I)$, and prove the so-called Sally-machine for the fiber cone.
Namely, assume that $a^0$ is filter-regular in $F(I)$ and
$a^{\ast}$ is filter-regular in $G(I)$. Then
$$\depth F(I/(a)) >0 \Rightarrow a^0 \text{ regular  in
} F(I).$$

Notice also that if in addition $a^{\ast} $ is regular in $G(I)$
then $F(I)/(a^0) \cong F(I/(a))$ and so $\depth F(I)= \depth
F(I/(a))+1$.

Now, we extend the Sally-machine for the fiber cone to sequences
of arbitrary length.

\begin{lem}
\label{SM1} Let $a_1,\dots ,a_k \in I$  such that $a_1^0,\dots
,a_k^0$ is a filter-regular sequence in $F(I)$ and $a_1^{\ast},
\dots ,a_{k}^{\ast}$ is a filter-regular sequence in $G(I)$.
Assume that $a_1^{\ast}, \dots ,a_{k-1}^{\ast}$ is regular in
$G(I)$. Then
$$ \depth F(I/(a_1,\dots ,a_k))>0
\Rightarrow a_1^0,\dots a_k^0 \text{ regular  in } F(I).$$ If in
addition $a_1^{\ast}, \dots ,a_{k}^{\ast}$ is a regular sequence
in $G(I)$ then $$\depth F(I)=depth F(I/(a_1,\dots ,a_k))+k\geq
1+k.$$
\end{lem}

\begin{proof}
The case $k=1$ is the above cited result as Sally-machine for the
fiber cone. Assume now $k>1$. For $i\leq k$ we will denote
$A_i=A/(a_1,\dots a_i)$, $\fm_i= \fm/(a_1,\dots ,a_i)$ and
$I_i=I/(a_1,\dots, a_i)$. Since $a_1^{\ast}, \dots
,a_{k-1}^{\ast}$ is regular in $G(I)$ we have that
$$G(I_{k-1})\simeq G(I)/(a_1^{\ast}, \ldots ,a_{k-1}^{\ast})$$
$$F(I_{k-1})\simeq F(I)/(a_1^0, \ldots ,a_{k-1}^0)$$ Therefore,
putting $\overline{a_k} \in A_{k-1}$ we have that
$\overline{a_k}^0 \in F(I_{k-1})$ and $\overline{a_k}^{\ast} \in
G(I_{k-1})$ are filter regular.

Assume now that $\depth F(I_k) >0$. Then, $F(I_k)\cong
F(I_{k-1}/(\overline{a_k}))$ and by induction for $k=1$ we get
that $\overline{a_k}^0$ is regular in $F(I_{k-1})$. Hence, $\depth
F(I_{k-1}) >0$ and again by induction for $k-1$, $a_1^0,\dots
,a_{k-1}^0$ is a regular sequence in $F(I)$, and $a_1^0,\dots
a_k^0 $ is a regular sequence in $F(I)$ as well.

For the last assertion, if we assume that $a_1^{\ast},\dots,
a_k^{\ast}$ is a regular sequence in $G(I)$, then $F(I_k)\cong
F(I)/(a_1^0,\dots ,a_k^0)$ with $a_1^0,\dots, a_k^0$ a regular
sequence in $F(I)$. So, $\depth F(I)=\depth F(I_k) +k\geq 1 + k$.
\end{proof}

Suppose now that $l:=l(I)\geq 1$ and let $J\subset I$ be a minimal
reduction of $I$. By \cite[Proposition 2.2]{JaVe2} there always
exist an element $a\in J\setminus \fm J$ such that $a^*$ is
filter-regular in $G(I)$ and $a^0$ is filter-regular in $F(I)$.
Moreover, if $\depth G_+(I)>0$, $a^*$ is a regular element in
$G(I)$ (\cite[Lemma 2.5]{Co}). Hence, if we assume that $\depth
G(I)_+ \geq l-1$, proceeding by induction one can always find
$(a_1, \ldots , a_l)=J$ such that $a_1^0, \ldots , a_l^0$ is a
filter-regular sequence in $F(I)$, $a_1^{\ast}, \ldots
,a_l^{\ast}$ is a filter regular sequence in $G(I)$, and
$a_1^{\ast}, \ldots ,a_{l-1}^{\ast}$ is a regular sequence in
$G(I)$. Set $A_{l-1}=A/(a_1, \ldots ,a_{l-1})$, $\fm_{l-1}=\fm
/(a_1, \ldots ,a_{l-1})$ and $I_{l-1}=I/(a_1, \ldots
,a_{l-1})\subset A_{l-1}$.

\begin{lem}
\label{SM11} $F(I)$ is Cohen-Macaulay if and only if $F(I_{l-1})$
is Cohen-Macaulay.
\end{lem}

\begin{proof}
First note that since $a_1^{\ast}, \ldots ,a_{l-1}^{\ast}$ is a
regular sequence in $G(I)$, then $$F(I)/(a_1^0, \ldots ,
a_{l-1}^0)\simeq F(I_{l-1})$$ and $\dim F(I_{l-1})=1$. Assume that
$F(I)$ is Cohen-Macaulay. Then, $a_1^0, \ldots , a_{l-1}^0$ is a
regular sequence in $F(I)$ and so $F(I_{l-1})$ is Cohen-Macaulay
too. Conversely, since $\depth F(I_{l-1})>0$ and $a_1^{\ast},
\ldots ,a_{l-1}^{\ast}$ is a regular sequence in $G(I)$ we have by
Lemma \ref{SM1} that $\depth F(I)= l-1+1=l$ and so $F(I)$ is
Cohen-Macaulay.
\end{proof}

Let $I$ be an ideal with $l(I)=l\geq 1$ and $\text{grade}(I)=l$,
that is, an equimultiple ideal with $\text{grade}(I)=l(I)$, and
assume that $\depth G(I)_+\geq l-1$. Then, the reduction number
$r_J(I)$ is independent of the choice of minimal reduction $J$
(see S. Huckaba \cite[Theorem 2.1]{Huc2}.

\begin{prop}
\label{SM2} (See \cite[Theorem 5.6 ] {HeKi}) Let $(A,\fm)$ be a
Noetherian local ring with an infinite residue field and let $I$
be an equimultiple ideal with analytic spread $l$,
$\text{grade}(I)= l$ and minimally generated by $l+1$ elements.
Assume that $\text{grade}(G(I)_+)\geq l-1$. Then $F(I)$ is
Gorenstein.
\end{prop}

\begin{proof}
Consider the ideal $I_{l-1}$ in the local ring $A_{l-1}$ defined
as above. This is a regular ideal of analytic spread one and
minimally generated by 2 elements. So, its fiber cone
$F(I_{l-1})\simeq F(I)/(a_1^0,\dots ,a_{l-1}^0)$ is Gorenstein by
Proposition \ref{D4}. Since $\depth F(I_{l-1})>0$, $a_1^0, \ldots
, a_{l-1}^0$ is a regular sequence by Lemma \ref{SM1} and $F(I)$
is also Gorenstein.
\end{proof}

From now on, given $a_1,\dots ,a_k \in A$ we will denote by
$A_k=A/(a_1,\dots ,a_k )$ and for an ideal $L$ of $A$ by
$L_k=(L+(a_1,\dots ,a_k ))/(a_1,\dots ,a_k )$.

Assume that $a_1,\dots ,a_k \in I$ is such that $a_1^{\ast},\dots
,a_k^{\ast} $ is a regular sequence in $G(I)$ and $a_1^0,\dots
,a_k^0$ is a regular sequence in $F(I)$. By the mixed
Valabrega-Valla criterium, see Cortadellas-Zarzuela \cite{CoZa1},
these conditions are equivalent to
\begin{enumerate}
\item $a_1,\dots , a_k$ is a regular sequence in $A$. \item
$(a_1,\dots ,a_{k})\cap I^n =(a_1,\dots ,a_{k})I^{n-1}$, for all
$n\geq 0$. \item $(a_1,\dots ,a_{k})\cap \fm I^n =(a_1,\dots
,a_{k})\fm I^{n-1}$ for all $n\geq 1$.
\end{enumerate}

\begin{lem} \label{E1}
Let $a_1,\dots ,a_k \in I$ such that $a_1^{\ast},\dots ,a_k^{\ast}
$ is a regular sequence in $G(I)$ and $a_1^0,\dots ,a_k^0$ is a
regular sequence in $F(I)$. Then,
$$ \mu (I_k^n)= \sum_{i=0}^k (-1)^i{k \choose i }\mu (I^{n-i}).$$
\end{lem}
\begin{proof} We use induction on $k$. For $k=1$, $\mu (I^n_1)= \la
(I^n/(\fm I^n +a_1I^n)$ by (2). Consider the exact sequence
$$0\rightarrow a_1I^{n-1}/(a_1I^{n-1}\cap \fm I^n) \rightarrow I^n/\fm
I^n \rightarrow I^n/(\fm I^n+a_1I^{n-1}) \rightarrow 0.$$ By
condition (3), $a_1I^{n-1}\cap \fm I^n=a_1\fm I^{n-1}$, and
$a_1I^{n-1}/a_1\fm I^{n-1} \cong I^{n-1}/\fm I^{n-1}$ by condition
(1). Thus $\mu (I_1^n)=\mu (I^n)- \mu (I^{n-1})$.

Let $1<k$. Then, $\overline a_k$, $(\overline a_k)^{\ast}$ and
$(\overline a_k)^0$ are regular elements, respectively, in the
rings $A_{k-1}$, $G(I_{k-1})$ and $F(I_{k-1})$. Therefore, $\mu
(I_k^n)=\mu (I_{k-1}^n)- \mu (I_{k-1}^{n-1})$ for the case $k=1$,
and by induction $$
\begin{array}{ll} \mu (I_k^n) & =\sum_{i=0}^{k-1} (-1)^i {k-1 \choose i } \mu
(I^{n-i})- \sum_{i=0}^{k-1} (-1)^i {k-1 \choose i } \mu
(I^{n-1-i})\\
  &= \sum_{i=0}^{k} (-1)^i ({k-1 \choose i }+{k-1 \choose i-1 } )\mu
 (I^{n-i})\\
&= \sum_{i=0}^k (-1)^i{k \choose i }\mu (I^{n-i}). \end{array}$$
\end{proof}

Let $(A,\fm)$ be a Cohen-Macaulay local ring of dimension $d$ and
$I$ be an $\fm$-primary ideal in $A$. By G. Valla \cite{Va},
$$
\la(I/I^2)=e(I)+(d-1)\la(A/I) -\la(I^2/JI)$$ for any $J$ minimal
reduction of $I$ where $e(I)$ denotes the multiplicity of the
ideal $I$. In particular, the length of $I^2/JI$ does not depend
of the minimal reduction of $I$.

\begin{prop}
\label{SM3}(See \cite{JaPuVe})  Let $(A,\fm)$ be a Cohen-Macaulay
local ring of dimension $d$ with an infinite residue field and let
$I$ be an $\fm$-primary ideal. Assume that $e(I)+(d-1)\la(A/I)
-\la(I/I^2)=1$. Let $J$ be a minimal reduction of $I$. Then $F(I)$
is a Cohen-Macaulay ring if and only if $\fm I^2 = J\fm I$. In
this case, $$\sum_{i=0}^d (-1)^i{d \choose i} \mu(I^{n+1-i})=1$$
for $1 \leq n \leq r(I)-1$.
\end{prop}

\begin{proof}
First, notice that the assumption $e(I)+(d-1)\la(A/I)
-\la(I/I^2)=1$ is equivalent to $\la (I^2/JI)=1$. By M. Rossi
\cite{Ro} $\depth G(I)\geq d-1$, hence there exist $(a_1,\dots
,a_d)=J$ such that $a_1^0,\dots,a_d^0$ is a filter regular
sequence in $F(I)$, $a_1^{\ast},\dots,a_d^{\ast}$ is a filter
regular sequence in $G(I)$, and $a_1^{\ast},\dots,a_{d-1}^{\ast}$
is a regular sequence in $G(I)$. By Lemma \ref{SM11} $F(I)$ is
Cohen-Macaulay if and only if $ F(I_{d-1})$ is Cohen-Macaulay.
Moreover, by S. Huckaba \cite[Lemma 1-1]{Huc2} $(\overline{a_d})$
is a minimal reduction of $I_{d-1}$ with $r(I_{d-1})=r(I)$, and
$\la (I_{d-1}^2/(\overline{a_d})I_{d-1})=\la(I^2/(JI+I^2\cap
(a_1,\dots,a_{d-1}))=\la(I^2/JI)=1$.

Suppose that $\fm I^2=J\fm I$. Then, $\fm_{d-1}
I_{d-1}^2=\overline{a_d}\fm_{d-1} I_{d-1}$ and $F(I_{d-1})$ is
Cohen-Macaulay by Proposition \ref{D5} and so $F(I)$ is
Cohen-Macaulay too. Conversely, if $F(I)$ is Cohen-Macaulay then
$F(I_{d-1})$ is also Cohen-Macaulay, and by Proposition \ref{D5}
$\fm_{d-1} I_{d-1}^2=\overline{a_d}\fm_{d-1} I_{d-1}$. Hence, $\fm
I^2 = J\fm I + \fm I^2\cap (a_1,\dots, a_{d-1}) = (a_1,\dots
,a_{d-1})\fm I=J\fm I$, since $a_1^{\ast},\dots,a_{d-1}^{\ast}$
and $a_1^0,\dots,a_{d-1}^0$ are regular sequences in $G(I)$ and
$F(I)$, respectively. Finally, if $F(I)$ is Cohen-Macaulay
$F(I_{d-1})$ is Cohen-Macaulay too and by Proposition \ref{D5}
$\mu(I_{d-1}^{n+1}) -\mu(I_{d-1}^n )=1$ for $1\leq n \leq r(I)-1$.
Therefore, we may apply Lemma \ref{E1} to obtain for $1\leq n \leq
r(I)-1$ the equality
$$1=\mu(I_{d-1}^{n+1}) -\mu(I_{d-1}^n )= \sum_{i=0}^d (-1)^i{d \choose i}
\mu(I^{n+1-i}).
$$
\end{proof}

\begin{thm}
\label{E2} Let $(A,\fm)$ be a Noetherian local ring with an
infinite residue field and let $I\subset A$ be an equimultiple
ideal with analytic spread $l$ and reduction number $r(I)$. Assume
that $\text{grade}(I)=l$, $\depth G(I) \geq l-1$ and $\depth
F(I)\geq l-1$. Let $J\subset I$ be a minimal reduction of $I$.
Then the following equalities hold:
\begin{enumerate}
\item $\mathrm{reg}(F(I))=r(I).$

\item $e(F(I))=\sum_{i=0}^{l-1} {l-1\choose i} \mu(I^{r-i})$.

\item $F(I)$ is a Cohen-Macaulay ring if and only if $$\la
(I^n/(\fm I^n+JI^{n-1})= \sum_{i=0}^{l} (-1)^i {l \choose i} \mu
(I^{n-i}),$$ for $1\leq n \leq r(I)$.
\end{enumerate}
\end{thm}

\begin{proof}
The results are true for $l=1$ by Theorem \ref{C5}. Assume that
$l\geq 1$ and let $(a_1, \ldots , a_l)=J$ such that $a_1, \ldots ,
a_l$ is a regular sequence in $A$, $a_1^{\ast}, \ldots ,
a_{l-1}^{\ast}$ is a regular sequence in $G(I)$, and $a_1^0,
\ldots , a_{l-}^0$ is a regular sequence in $F(I)$. Observe that
$I_{l-1}\subset A_{l-1}$ is a regular ideal with analytic spread
$1$ and the same reduction number as $I$. Then,
$\mathrm{reg}(F(I)) =\mathrm{reg}(F(I)/(a_1^0,\dots ,a_{l-1}^0)) =
\mathrm{reg}(F(I_{l-1})) =r(I)$. Similarly,
$e(F(I))=e(F(I_{l-1}))= \mu (I_{l-1}^{r})=\sum_{i=0}^{l-1}(-1)^i {
l-1 \choose i} \mu (I^{r-i})$ by Lemma \ref{E1}.

Now, by Lemma \ref{SM11} $F(I)$ is  Cohen-Macaulay  if, and only
if, $F(I_{l-1})$ is Cohen-Macaulay, and by Theorem \ref{C5} this
is equivalent to $\la (I_{l-1}^n/(\fm_{l-1} I_{l-1}^n
+a_lI_{l-1}^{n-1}))=\mu (I_{l-1}^n)-\mu (I_{l-1}^{n-1})$, for all
$1\leq n\leq r(I)$. Then, we get (3) by the isomorphisms
$$ \begin{array}{ll}
I_{l-1}^n/(\fm I_{l-1}^n +a_lI_{l-1}^{n-1}) & \cong I^n/(\fm I^n
+a_lI^{n-1} +I^n \cap(a_1, \dots, a_{l-1})) \\
&\cong I^n/(\fm I^n + (a_1, \dots, a_l))I^{n-1}), \end{array}$$
and taking into account that $\mu (I_{l-1}^n)-\mu (I_{l-1}^{n-1})=
\sum_{i=0}^{l} (-1)^i {l \choose i} \mu (I^{n-i})$ by Lemma
\ref{E1}.
\end{proof}

\begin{rem}
\label{E3} It may be seen that if $I$ is an equimultiple ideal
with analytic spread $l(I)$ and reduction number $r(I)$, such that
$\text{grade}(I)=l(I)$ and $\depth G(I) \geq l(I)-1$, then
$\mathrm{reg}(G(I))=r(I)$, see Hoa-Zarzuela \cite[Proposition
3.6]{HoZa}.
\end{rem}

\bibliographystyle{amsalpha}

 \end{document}